\documentclass{amsart}

\usepackage{doc}
\usepackage{amssymb}
\usepackage{amsmath}
\usepackage{amsthm}

\newtheorem{theorem}{Theorem}
\newtheorem{corollary}{Corollary}
\newtheorem{proposition}{Proposition}
\theoremstyle{definition}
\newtheorem{definition}{Definition}

\newcommand{\SSigma}[1]{\boldsymbol{\Sigma_{#1}}}
\newcommand{\PPi}[1]{\boldsymbol{\Pi_{#1}}}
\newcommand{\enf}[1]{\left\|{#1}\right\|}

\title{Applications of the analogy between formulas and exponential
  polynomials to equivalence and normal forms}

\author{Danko Ilik}

\begin{document}

\maketitle

\begin{abstract}
  We show some applications of the formulas-as-polynomials
  correspondence: 1) a method for (dis)proving formula isomorphism and
  equivalence based on showing (in)equality; 2) a constructive
  analogue of the arithmetical hierarchy, based on the exp-log normal
  form. The results are valid intuitionistically, as well as
  classically.
\end{abstract}

\section{Introduction}
\label{sec:intro}

The language of propositional formulas bears a striking similarity to
the languages of exponential polynomials, when one identifies
$ \phi\vee \psi$ with $\phi+\psi$, $\phi\wedge \psi$ with
$\phi\times \psi$, $\psi\Rightarrow \phi$ with $\phi^{\psi}$, and
atomic propositions $\chi_i$ with variables. This algebraic notation
for formulas is often used in categorical logic and type
theory~\cite{lambek2011reflections}, however some of the deeper
logical implications of the analogy between formulas and exponential
polynomials are seldom pointed out. In this paper, we review some of
the logical implications in relation to formula equivalence, formula
isomorphism, and formula normal forms.

We first deal with the propositional case and later, in
Section~\ref{sec:hierarchy}, we will extend the analogy between
formulas and exponential polynomials to the first-order case. We work
in the context of intuitionistic logic, but, since the results are
about formula equivalence, they are also relevant classically.

Central will be the notion of formula isomorphism, defined as an
extension of intuitionistic equivalence---purely logical equivalence
with no use of extra-logical axioms such as induction, choice,
Church's thesis, etc.---with a property expressing that the equivalent
formulas, when seen as sets of their proofs, have isomorphic
structures. To be concrete, we will use the language of typed lambda
calculus and the associated $\beta$-$\eta$-equality relation to denote
extensional equivalence between natural deduction trees that represent
proofs. We define the isomorphism between formulas $\phi$ and $\psi$
by
\begin{align*}
  \phi\cong\psi && \text{iff} && \text{ there exist proofs } M : (\Gamma\vdash\phi\Rightarrow\psi) \text{ and } N : (\Gamma\vdash\psi\Rightarrow\phi)\text{ such that }\\
                && &&\lambda x. M (N x) =_{\beta\eta} \lambda x. x \text{ and } \lambda y. N (M y) =_{\beta\eta} \lambda y. y.
\end{align*}
This definition of isomorphism via typed lambda terms is equivalent to
the categorical approach to proof representations; see for instance
\cite{Dosen2003}. What is meant is that, not only does $M$ prove the
implication $\phi\Rightarrow\psi$ and $N$ prove $\psi\Rightarrow\phi$,
but also that any proof $x$ of $\psi$ can be mapped to a proof $N x$
of $\phi$ and then \emph{back} to the \emph{same} proof $x$ of $\psi$, that
is, without loss of information; and analogously for any proof $y$ of
$\phi$. For a simple example of formulas that are equivalent but not
isomorphic, consider the equivalence
$\alpha\wedge \alpha\Leftrightarrow \alpha$, where $\alpha$ has for
instance exactly 2 possible proofs.

Some rationale for why formula isomorphism is the ``right'' form of
equivalence for intuitionistic logic is given in
\cite{Ilik:2019:hierarchy}, but in this paper we are primarily
interested in the notion because of its link to exponential
polynomials, explained in the next section.


\section{(Dis)proving propositional formula isomorphism}
\label{sec:disproving}

We start with a well known fact from categorical logic and typed
lambda calculi (e.g.~\cite{FioreDCB2006}),
\[
  \phi\cong\psi\quad\to\quad \mathbb{N}^+\vDash \phi=\psi,
\]
linking isomorphism of formulas to realizability in~$\mathbb{N}^+$.
The expression on the right is the standard model theoretic one: for
any instantiation of the variables with positive natural numbers,
$\phi$ and $\psi$ compute to the same number\footnote{This number is
  always positive, because the exponential polynomials that we need
  for the analogy to formulas have positive coefficients only.}. This
immediately gives a method for disproving formula isomorphism.
\begin{corollary}
  To show that $\phi\not\cong\psi$, it is enough to show (by whatever
  means, i.e., algebra, analysis, etc.) that $\phi\neq\psi$.
\end{corollary}

Although the method directly follows from basic relation between
isomorphism and equality, we think that it is interesting to mention,
because it can be much simpler to employ than:
\begin{itemize}
\item showing that there is no proof term for $\phi\Rightarrow\psi$ or
  a proof term for $\psi\Rightarrow\phi$;
\item showing that no proof term for $\phi\Rightarrow\psi$ and proof
  term for $\psi\Rightarrow\phi$ can be $\beta$-$\eta$-equal;
\item proving that $\phi\Leftrightarrow\neg\psi$ (or vice versa); or
\item building a model that realizes $\phi$ but does not realize
  $\psi$ (or vice versa).
\end{itemize}

Of course, it is not a method for disproving \emph{equivalence},
because even when $\phi\not\cong\psi$, $\phi\Leftrightarrow\psi$ may
still hold. But, the following link between equality and isomorphism
allows to \emph{prove} equivalence:
\[
  \text{HSI}\vdash \phi ~\dot=~ \psi\quad\to\quad \phi\cong\psi.
\]
The expression $\text{HSI}\vdash\phi~\dot=~\psi$ means that the
equality between the exponential polynomials $\phi$ and $\psi$ is
\emph{derivable} from the so called high-school identities:
\begin{align*}
  \phi &= \phi              & \phi 1 &= \phi\\
  \phi+\psi &= \psi+\phi          & \phi^1&= \phi\\
  (\phi+\psi)+\xi &= \phi+(\psi+\xi)  & 1^\phi&= 1\\
  \phi \psi &= \psi \phi          & \phi^{\psi+\xi}&= \phi^\psi \phi^\xi\\
  (\phi \psi) \xi &= \phi (\psi \xi)  & (\phi \psi)^\xi&= \phi^\xi \psi^\xi\\
  \phi(\psi+\xi) &= \phi \psi + \phi \xi & (\phi^\psi)^\xi&= \phi^{\psi \xi}.
\end{align*}
\begin{corollary}
  To show that $\phi\cong\psi$, it is enough to show that $\phi=\psi$
  is derivable by the high-school identities.
\end{corollary}

Note that not every true equality between exponential polynomials is
derivable by the high-school identities (HSI),
\[
  \mathbb{N}^+\vDash\phi=\psi\quad\not\to\quad\text{HSI}\vdash \phi~\dot=~\psi.
\]
This is a consequence of the negative solution to Tarski's high-school
algebra problem~\cite{burris04}. A typical operation that one can
perform in the model, but not in the derivation, is canceling.  Here is
a simple example of a $\phi$ and a $\psi$, such that
$\mathbb{N}^+\vDash\phi=\psi$ and $\phi\cong\psi$, but
$\text{HSI}\not\vdash \phi~\dot=~\psi$, due to Martin~\cite{martin73}:
\begin{align*}
  \left(\chi_4\Rightarrow\left(\chi_3\Rightarrow \chi_1\right)\vee\left(\chi_3\Rightarrow \chi_1\right)\right)\wedge\left(\chi_3\Rightarrow\left(\chi_4\Rightarrow \chi_2\right)\vee\left(\chi_4\Rightarrow \chi_2\right)\right),\\
  \left(\chi_3\Rightarrow\left(\chi_4\Rightarrow \chi_1\right)\vee\left(\chi_4\Rightarrow \chi_1\right)\right)\wedge\left(\chi_4\Rightarrow\left(\chi_3\Rightarrow \chi_2\right)\vee\left(\chi_3\Rightarrow \chi_2\right)\right).
\end{align*}

It is possible to ``repair'' HSI by extending it to an enumerable set
of axioms HSI* due to Wilkie~\cite{wilkie00}, that, in addition to the
HSI, contains all true equalities between \emph{ordinary} (i.e.,
non-exponential) positive polynomials, possibly with negative
coefficients. Then, we do get that
\begin{align*}
  \mathbb{N}^+\vDash\phi=\psi&&\leftrightarrow&&\quad\text{HSI}^*\vdash \phi~\dot=~\psi\\
  \phi\cong\psi&&\rightarrow&&\quad\text{HSI}^*\vdash \phi~\dot=~\psi
\end{align*}
however it is not clear how to prove in general that
$\text{HSI}^*\vdash \phi~\dot=~\psi$ implies $\phi\cong\psi$, since
the ordinary polynomial equalities we have added to HSI* may contain
negative coefficients. A partial solution is given
in~\cite{ilik:2014:axioms}.

Finally, we remark that, although the analogy between isomorphism and
equality can be delicate for certain kind of formulas---the formulas
of the kind of Martin---very often it works seamlessly. There are also
some general results in this direction based on the form of formulas
in question, 
such as
the following one.
\begin{proposition}[\cite{gurevic93,levitz75,ilik:2014:axioms}]
  For the class $\mathcal{L}$ of formulas of Gurevi\v{c}--Levitz, we
  have that, for all $\phi,\psi\in\mathcal{L}$,
  \[
    \text{HSI}\vdash \phi~\dot=~\psi ~\leftrightarrow~ \phi\cong\psi ~\leftrightarrow~ \mathbb{N}^+\vDash\phi=\psi,
  \]
  where the class $\mathcal{L}$ is defined inductively:
  \begin{align*}
    \mathcal{L}\ni\phi,\psi &::= \chi_i ~|~ \phi\vee\psi ~|~ \phi\wedge\psi ~|~ \phi\Rightarrow\lambda\\
    \Lambda\ni\lambda,\mu,\nu &::= \pi_i ~|~ \mu\vee\nu ~|~ \mu\wedge\nu ~|~ \mu\Rightarrow \lambda_0,
  \end{align*}
  where $\lambda_0\in\Lambda$ but $\lambda_0$ contains no proposition
  (i.e., variable) $\pi_i$.
\end{proposition}
For these kinds of formulas, equivalence and isomorphism can be proved
by showing equality of exponential polynomials (by whichever
mathematical means) and not only by deriving equality via the axioms
HSI.

\section{Intuitionistic ``prenex'' normal form}
\label{sec:hierarchy}

For ordinary (i.e., non-exponential) polynomials, equality is
decidable and so is isomorphism. This fragment of polynomials
corresponds to the fragment of formulas constructed from the two
logical connectives $\{\vee,\wedge\}$. If we consider the fragment
$\{\Rightarrow,\wedge\}$, the situation is the same. In both cases, it
is a consequence of the fact that there is a \emph{canonical} normal
form for the formulas (polynomials) in question.

In the general case when all of the connectives
$\{\wedge,\vee,\Rightarrow\}$ are present, a canonical normal form is
not known. However, we sometimes find~\cite{hardy} the following
transformation
\[
  \sigma^\tau = \exp{(\tau \log{\sigma})}
\]
that allows to derive a quasi\footnote{The normal form is ``quasi'',
  because, unlike the normal form for ordinary polynomials, it is not
  canonical, i.e., there are equal exponential polynomials that do not
  have the same normal form.} normal form for exponential polynomials.
\begin{theorem}[Theorem~2.1 of \cite{Ilik:2019:hierarchy}]
  Every propositional formula $\phi$ can be normalized to a formula
  $\enf{\phi}$, such that $\phi\cong\enf{\phi}$ and
  $\enf{\phi}\in\mathbf{\Pi}\cup\mathbf{\Sigma}$, where the classes
  $\mathbf{\Pi}$ and $\mathbf{\Sigma}$ are defined inductively and
  mutually as follows:
  \begin{align*}
    \mathbf{\Pi}\ni \gamma &::= (\gamma_1\Rightarrow\beta_1)\wedge\cdots\wedge(\gamma_n\Rightarrow\beta_n) & (n\ge 0)~\\
    \mathbf{\Sigma}\ni \beta &::= \chi_i ~|~ \gamma_1\vee\cdots\vee\gamma_n & (n\ge 2),
  \end{align*}
  where $\chi_i$ are prime formulas.
\end{theorem}

To extend the normal form to the first-order quantifiers, we adopt an
\emph{extended} exponential polynomial notation: we write
$\exists x \phi$ as $x \phi$ and $\forall x \phi$ as $\phi^x$, the
distinction between conjunctions and existential quantifiers, and
implications and universal quantifiers, being made by a variable
convention: we ``left-multiply'' and ``exponentiate'' by the lowercase
Latin letters $x, y, z$ in order to express quantifiers, while if we
do it with Greek $\phi, \psi$, it means that we are making a
conjunction and implication. Using this notation, the following
formula isomorphisms acquire the form of equations:
\begin{align*}
  \forall x (\phi\wedge \psi)&\cong \forall x \phi\wedge\forall x \psi &(\phi \psi)^x &= \phi^x \psi^x\\
  \exists x (\phi\vee \psi)&\cong \exists x \phi \vee \exists x \psi & x(\phi+\psi) &= x\phi + x\psi\\
  \exists x \phi\Rightarrow \psi&\cong \forall x(\phi\Rightarrow \psi) & \psi^{x\phi} &= (\psi^\phi)^x & (\text{where } x\notin\text{FV}(\psi))\\
  \psi\Rightarrow \forall x \phi&\cong \forall x(\psi\Rightarrow \phi) & (\phi^x)^\psi &= (\phi^\psi)^x & (\text{where } x\notin\text{FV}(\psi))
\end{align*}
If we add these equations to the axioms of HSI, we will thus still
preserve isomorphism. However, we do not know whether a relation to
model theoretic equality will hold, such as the one in the
propositional case. The reason is that we have extended the language
of exponential polynomials, and one needs to find an interpretation for
the new arithmetic operations $x(\cdot)$ and $(\cdot)^x$.
Nevertheless, the representation of formulas as extended exponential
polynomials is sufficient to show a normal form theorem for
first-order formulas.
\begin{theorem}[Theorem~4.1 of \cite{Ilik:2019:hierarchy}]
  Every first-order formula $\phi$ can be normalized to a formula
  $\enf{\phi}$, such that $\phi\cong\enf{\phi}$ and
  $\enf{\phi}\in\mathbf{\Pi}\cup\mathbf{\Sigma}$, where the classes
  $\mathbf{\Pi}$ and $\mathbf{\Sigma}$ are defined inductively and
  mutually as follows:
  \begin{align*}
    \mathbf{\Pi}\ni \gamma &::= \forall x_1(\gamma_1\Rightarrow\beta_1)\wedge\cdots\wedge\forall x_n(\gamma_n\Rightarrow\beta_n) & (n\ge 0)~\\
    \mathbf{\Sigma}\ni \beta &::= \chi_i ~|~ \gamma_1\vee\cdots\vee\gamma_n ~|~ \exists x \gamma& (n\ge 2),
  \end{align*}
  where $\chi_i$ are prime formulas.
\end{theorem}

Intuitionistically, this theorem is interesting, because we do not
have a notion of arithmetical hierarchy as versatile as the one of
classical logic; see~\cite{Ilik:2019:hierarchy} for a
discussion. Classically, one can also use this hierarchy as an
alternative to the standard one.

To make the link to the classical arithmetical hierarchy, we will
first need to assign levels to the hierarchy from the last theorem.
\begin{definition}
  The \emph{intuitionistic arithmetical hierarchy} is defined by
  assigning levels, $\SSigma{n}, \PPi{n}$, for $n\in\mathbb{N}$, to
  the formula classes $\SSigma{}$ and $\PPi{}$, in the following way:
  \begin{align*}
    \PPi{0}\ni \gamma &::= \top\Rightarrow\chi & \chi\text{ is a prime formula }\\
    \SSigma{0}\ni \beta &::= \chi & \chi\text{ is a prime formula }\\
    \PPi{n+1}\ni \gamma &::= \forall x_1(\gamma_1\Rightarrow\beta_1)\wedge\cdots\wedge\forall x_m(\gamma_m\Rightarrow\beta_m) & n= \max_{i=1}^m\{k ~|~ \beta_i\in\SSigma{k}\}\\
    \SSigma{n+1}\ni \beta &::= \chi_i ~|~ \gamma_1\vee\cdots\vee\gamma_m ~|~ \exists x \gamma & n= \max_{i=1}^m\{k ~|~ \gamma_i\in\PPi{k}\} \text{ or } \gamma\in\PPi{n}.
  \end{align*}
  We also extend the relation ``$\in$'' from formulas satisfying the
  inductive definition to all formulas, in the following way:
  $F\in\PPi{n}$ iff $\enf{F}\in\PPi{n}$; $F\in\SSigma{n+1}$ iff
  $\enf{F}\in\SSigma{n+1}$.
\end{definition}

The relation to the classical arithmetical hierarchy is given by the
following theorem, where the levels of the classical hierarchy are
denoted $\Sigma^0_n$ and $\Pi^0_n$, that is, with a zero superscript
and no bold face.
\begin{theorem}[Propositions 4.7 and 4.8 of \cite{Ilik:2019:hierarchy}]
  Say that a formula $\phi$ is \emph{classically represented in
    $\SSigma{n}$ (or $\PPi{n}$)} when there is a formula
  $\phi'\in\SSigma{n}$ (or $\PPi{n}$) such that $\phi$ and $\phi'$ are
  classically equivalent.

  If $\phi\in\Sigma^0_n$, then $\phi$ is classically represented in
  $\SSigma{n}$. If $\phi\in\Pi^0_n$, then $\phi$ is classically
  represented in $\PPi{n}$.

  Suppose that $\psi$ is in prenex normal form and with alternating
  quantifiers. Then: $\psi\in\SSigma{n}$ implies $\psi\in\Sigma^0_n$;
  $\psi\in\PPi{n}$ implies $\psi\in\Pi^0_n$.
\end{theorem}

As a corollary, we get that the intuitionistic hierarchy is proper.
\begin{corollary}[Corollary 4.9 of \cite{Ilik:2019:hierarchy}]
  For $n\ge 0$, $\SSigma{n}\subsetneq\SSigma{n+1}$,
  $\SSigma{n}\subsetneq\PPi{n+1}$, $\PPi{n}\subsetneq\SSigma{n+1}$,
  and $\PPi{n}\subsetneq\PPi{n+1}$.
\end{corollary}

One can think of the exp-log normal form of formulas as the
constructive analogue of the prenex normal form from classical logic.




\label{sec:bib}
\bibliographystyle{plain}
\bibliography{formpoly}
\nocite{*}

\end{document}